\documentclass{birkjour}
 \newtheorem{thm}{Theorem}[section]
 \newtheorem{cor}[thm]{Corollary}
 \newtheorem{lem}[thm]{Lemma}
 
 \theoremstyle{definition}
 \newtheorem{defn}[thm]{Definition}
 \theoremstyle{remark}
 \newtheorem{rem}[thm]{Remark}
 \newtheorem*{ex}{Example}
 \numberwithin{equation}{section}

\begin{document}

%
%
%
%
%
%
%
%
%

\title[Trace cohomology revisited]{Trace cohomology revisited}

\author[Nikolaev]{Igor V. Nikolaev}

\address{%
Department of Mathematics and Computer Science\\
St.~John's University, 8000 Utopia Parkway\\
New York,  NY 11439,  United States}

\email{igor.v.nikolaev@gmail.com}




\date{January 1, 2004}

\begin{abstract}
We use a cohomology theory coming from the canonical trace
on a $C^*$-algebra of the projective variety to prove an analog 
of the Riemann Hypothesis for the Kuga-Sato varieties. 
\end{abstract}

\maketitle

\section{Introduction}
The aim of our note is a proof of the Riemann Hypothesis for a class of 
projective varieties  over finite fields using the notion of a trace cohomology
introduced in \cite{Nik1}.    To define such a cohomology, recall that 
the  Serre $C^*$-algebra  $\mathcal{A}_V$ of an $n$-dimensional complex projective variety
$V_\mathbf{C}$ is the norm-closure of a self-adjoint representation of the 
twisted homogeneous coordinate ring of $V_\mathbf{C}$  by the bounded 
linear operators on a Hilbert space $\mathcal{H}$.  
We shall write  $\tau: \mathcal{A}_V\otimes \mathcal{K}\to   \mathbf{R}$
to denote   the canonical  normalized trace on  the stable $C^*$-algebra $\mathcal{A}_V\otimes \mathcal{K}$,  
 i.e. a positive linear functional   of norm $1$,   such that $\tau(yx)=\tau(xy)$ for all $x,y\in \mathcal{A}_V\otimes \mathcal{K}$,  see  
 [Blackadar 1986] \cite{B},  p. 31. 
Applying the Chern character formula to the  algebra $\mathcal{A}_V\otimes \mathcal{K}$,  
one obtains an injective homomorphism  $\tau_*: H^i(V_\mathbf{C})\longrightarrow \mathbf{R}$.
The $i$-th   trace cohomology group $\{H^i(V_\mathbf{C}) ~|~  0\le i\le 2n\}$ of  $V_{\mathbf{C}}$ is an
additive abelian subgroup $\tau_*(H^i(V_\mathbf{C}))$ of the real line $\mathbf{R}$. 
We refer the reader to Section 2 for the details.

It was shown in \cite{Nik2}  that   all of Weil's Conjectures, except for  an analog of the Riemann Hypothesis (RH), 
follow from simple properties of the trace cohomology.  
Recall that the Kuga-Sato variety is a fiber product of the modular curves;
we refer the reader to  Section 2.4 for an exact definition. 
The $i$-th cohomology group of such a variety is related to the space of cusp forms of weight $i+1$
[Deligne 1969]  \cite{Del1}  and [Scholl 1985]  \cite{Sch1}.
In this note we use the trace cohomology to prove the RH
for the Kuga-Sato varieties with  a lifting from characteristic $p$ to characteristic zero 
 [Hartshorne 2010]  \cite[Theorem 22.1]{H2}.  Namely,
 let ${\Bbb F}_q$ be a finite field with $q=p^r$ elements and $V({\Bbb F}_q)$ be a
smooth $n$-dimensional Kuga-Sato variety over ${\Bbb F}_q$. 
\begin{thm}\label{thm1}
The roots $\alpha_{ij}$ of polynomials $P_i(t)$ in the zeta 
function $Z_V(t)={P_1(t)\dots P_{2n-1}(t)\over P_0(t)\dots P_{2n}(t)}$
of $V({\Bbb F}_q)$
are algebraic numbers of the absolute value  $|\alpha_{ij}|=q^{{i\over 2}}$. 
\end{thm}

Theorem \ref{thm1} is not new. It has been proved in full generality 
and different methods in the classical work   [Deligne 1974] \cite{Del2}.
The novelty of our approach are concepts of noncommutative 
geometry,  e.g.  the Serre $C^*$-algebras and trace cohomology. 
The latter provide pathways and tools to some  open problems 
of algebraic geometry and number theory inaccessible otherwise. 

\medskip
The article is organized as follows.  The Serre $C^*$-algebras and trace
cohomology are introduced in Section 2.  Theorem  \ref{thm1} 
is  proved in Section 3.  The trace cohomology of an algebraic
 curve  is  calculated in Section 4.

\section{Preliminaries}
An excellent  survey of noncommutative algebraic geometry is  written by  [Stafford \& van ~den ~Bergh 2001]  \cite{StaVdb1}. 
For an introduction to  the $C^*$-algebras and their $K$-theory  we refer the reader to   [Murphy 1990]  \cite{M}
and  [Blackadar 1986]  \cite{B},  respectively. 
The Serre $C^*$-algebras were defined  in \cite{Nik1} and the  trace cohomology in \cite{Nik2}.  
The Weil's Conjectures were introduced in [Weil 1949]
\cite{Wei1}.

\subsection{Serre $C^*$-algebra}
Let $V$ be a projective scheme over a field $k$ and let $\mathcal{L}$ be an invertible
sheaf of linear forms on $V$.  If $\sigma$ is an automorphism of $V$,  then
the pullback of $\mathcal{L}$ along $\sigma$ will be denoted by $\mathcal{L}^{\sigma}$,
i.e. $\mathcal{L}^{\sigma}(U):= \mathcal{L}(\sigma U)$ for every $U\subset V$. 
Consider the graded $k$-algebra
\begin{equation}
B(V, \mathcal{L}, \sigma)=\bigoplus_{i\ge 0} H^0(V, ~\mathcal{L}\otimes \mathcal{L}^{\sigma}\otimes\dots
\otimes  \mathcal{L}^{\sigma^{ i-1}})
\end{equation}
called a {\it twisted homogeneous coordinate ring} of $V$;  notice that such a ring is 
non-commutative unless $\sigma$ is the trivial automorphism.    Recall that 
multiplication of sections of  $B(V, \mathcal{L}, \sigma)$ is defined by the 
rule  $ab=a\otimes b^{\sigma^m}$,   where $a\in B_m$ and $b\in B_n$.
Given a pair $(V,\sigma)$ consisting of a Noetherian scheme $V$ and 
 an automorphism $\sigma$ of $V$,  an invertible sheaf $\mathcal{L}$ on $V$
 is called {\it $\sigma$-ample}, if for every coherent sheaf $\mathcal{F}$ on $V$,
 the cohomology group $H^q(V, ~\mathcal{L}\otimes \mathcal{L}^{\sigma}\otimes\dots
\otimes  \mathcal{L}^{\sigma^{ n-1}}\otimes \mathcal{F})$  vanishes for $q>0$ and
$n>>0$.  Notice,  that if $\sigma$ is trivial,  this definition is equivalent to the
usual definition of ample invertible sheaf  [Serre 1955]  \cite{Ser1}.    
If $\mathcal{L}$ is a $\sigma$-ample invertible sheaf on $V$,  then
\begin{equation}\label{eq5}
Mod~(B(V, \mathcal{L}, \sigma)) ~/~Tors \cong Coh~(V),
\end{equation}
where  $Mod$ is the category of graded left modules over the ring $B(V, \mathcal{L}, \sigma)$,
$Tors$ is the full subcategory of $Mod$ of the torsion
modules and  $Coh$ is the category of quasi-coherent sheaves on a scheme $V$,
see [M.~Artin \& van den Bergh  1990]  \cite{ArtVdb1}.  In view of (\ref{eq5}) 
the ring $B(V, \mathcal{L}, \sigma)$  is indeed a coordinate ring of $V$,   see   [Serre 1955]  \cite{Ser1}.
\begin{rem}
Suppose that $R$ is  a commutative  graded ring,   such that $V=Spec~(R)$  is a
projective variety. 
Denote by $R[t,t^{-1}; \sigma]$  the ring of skew Laurent polynomials defined by the commutation 
relation  $b^{\sigma}t=tb$  for all $b\in R$,  where $b^{\sigma}$ is the image of $b$ under automorphism 
$\sigma: V\to V$.    Then $R[t,t^{-1}; \sigma]\cong B(V, \mathcal{L}, \sigma)$,  
see [M.~Artin \& van den Bergh  1990]  \cite{ArtVdb1}.
\end{rem}
Let $\mathcal{H}$ be a Hilbert space and   $\mathcal{B}(\mathcal{H})$ the algebra of 
all  bounded linear  operators on  $\mathcal{H}$.  For a  ring of skew Laurent polynomials $R[t, t^{-1};  \sigma]$,  
we shall consider a homomorphism 
\begin{equation}\label{eq6}
\rho: R[t, t^{-1};  \sigma]\longrightarrow \mathcal{B}(\mathcal{H}). 
\end{equation}
Recall  that the algebra $\mathcal{B}(\mathcal{H})$ is endowed  with a $\ast$-involution
 coming  from the scalar product on the Hilbert space $\mathcal{H}$. 
We shall call representation (\ref{eq6})  {\it $\ast$-coherent}  if
(i)  $\rho(t)$ and $\rho(t^{-1})$ are unitary operators,  such that
$\rho^*(t)=\rho(t^{-1})$ and   (ii) for all $b\in R$ it holds $(\rho^*(b))^{\sigma(\rho)}=\rho^*(b^{\sigma})$, 
where $\sigma(\rho)$ is an automorphism of  $\rho(R)$  induced by $\sigma$. 
Whenever  $B:=R[t, t^{-1};  \sigma]$  admits a $\ast$-coherent representation,
$\rho(B)$ is a $\ast$-algebra;  the norm-closure of  $\rho(B)$  yields
a   $C^*$-algebra,  see e.g.  [Murphy 1990]  \cite{M},  Section 2.1.  
\begin{defn}
 By  a  Serre $C^*$-algebra  of $V$ we understand the norm-closure of  $\rho(B)$;
 such a $C^*$-algebra   will be denoted  by  $\mathcal{A}_V$. 
\end{defn}
\begin{rem}\label{rmk2}
Each Serre $C^*$-algebra $\mathcal{A}_V$ is a crossed product $C^*$-algebra, 
see e.g.  [Williams  2007]  \cite{W},  pp 47-54 for the definition and details; 
namely,   $\mathcal{A}_V\cong C(V)\times_{\sigma} \mathbf{Z}$,  where $C(V)$ is the
$C^*$-algebra of all continuous complex-valued functions on $V$
and $\sigma$  is a $\ast$-coherent  automorphism of  $V$.   
\end{rem}
\begin{rem}
Let $\mathcal{K}$ be the $C^*$-algebra of all compact operators on a Hilbert space $\mathcal{H}$.
The stable Serre $C^*$-algebra $\mathcal{A}_V\otimes \mathcal{K}$ is endowed  with the unique normalized
trace (tracial state) $\tau:   \mathcal{A}_V\otimes \mathcal{K}\to \mathbf{R}$,   i.e. a positive linear functional
of norm $1$  such that $\tau(yx)=\tau(xy)$ for all $x,y\in \mathcal{A}_V\otimes \mathcal{K}$,  see  
 [Blackadar 1986] \cite{B},  p. 31. 
\end{rem}

\subsection{Trace cohomology}
Let $k$ be a number field.   Let $V(k)$ be  a smooth $n$-dimensional projective 
variety over $k$,  such  that  variety $V:=V({\Bbb F}_q)$ is the reduction modulo $q$ of $V(k)$
for a fixed choice of integral model [Hartshorne 2010,  Theorem 22.1]  \cite{H2}.
In other words,  $V(k)$ is defined by polynomial equations for $V$  over the field of 
complex numbers.  
Because the Serre $C^*$-algebra $\mathcal{A}_V$ of $V(k)$ is a crossed product $C^*$-algebra of the form
$\mathcal{A}_V\cong C(V(k))\times \mathbf{Z}$ (Remark \ref{rmk2}),  one can use  the Pimsner-Voiculescu 
six term exact sequence for the crossed products, see  e.g.  [Blackadar 1986]  \cite{B}, p. 83 for
 the details.  Thus   one gets the  short exact sequence of the algebraic $K$-groups:  
$0\to K_0(C(V(k)))\buildrel  i_*\over\to  K_0(\mathcal{A}_V)\to K_1(C(V(k)))\to 0$, 
where   map  $i_*$  is induced by the natural  embedding of $C(V(k))$ 
into $\mathcal{A}_V$.   We  have $K_0(C(V(k)))\cong K^0(V(k))$ and 
$K_1(C(V(k)))\cong K^{-1}(V(k))$,  where $K^0$ and $K^{-1}$  are  the topological
$K$-groups of the variety $V(k)$, see  [Blackadar 1986]  \cite{B}, p. 80. 
By  the Chern character formula,  one gets
$K^0(V(k))\otimes \mathbf{Q}\cong H^{even}(V(k); \mathbf{Q})$ and 
$K^{-1}(V(k))\otimes \mathbf{Q}\cong H^{odd}(V(k); \mathbf{Q})$, 
where $H^{even}$  ($H^{odd}$)  is the direct sum of even (odd, resp.) 
cohomology groups of $V(k)$.  
Notice that $K_0(\mathcal{A}_V\otimes \mathcal{K})\cong K_0(\mathcal{A}_V)$ because
of  stability of the $K_0$-group with respect to tensor products by the algebra 
$\mathcal{K}$,  see e.g.   [Blackadar 1986]  \cite{B}, p. 32.
One gets the   commutative diagram in Figure 1, 
where $\tau_*$ denotes  a homomorphism  induced on $K_0$ by  the canonical  trace 
$\tau$ on the $C^*$-algebra  $\mathcal{A}_V\otimes \mathcal{K}$.  
Recall that   $H^{even}(V(k)):=\oplus_{i=0}^n H^{2i}(V(k))$ and  
$H^{odd}(V):=\oplus_{i=1}^n H^{2i-1}(V(k))$, where $H^*(V(k))$ is the
singular cohomology of $V(k)$. 
 One gets  for each  $0\le i\le 2n$ 
 an injective  homomorphism 
 \begin{equation}\label{eq7}
 \tau_*:  ~H^i(V(k))\longrightarrow  \mathbf{R}. 
 \end{equation}
\begin{figure}
\begin{picture}(300,100)(0,0)
\put(160,72){\vector(0,-1){35}}
\put(80,65){\vector(2,-1){45}}
\put(240,65){\vector(-2,-1){45}}
\put(10,80){$ H^{even}(V(k))\otimes \mathbf{Q} 
\buildrel  i_*\over\longrightarrow  K_0(\mathcal{A}_V\otimes\mathcal{K})\otimes \mathbf{Q} 
\longrightarrow H^{odd}(V(k))\otimes \mathbf{Q}$}
\put(167,55){$\tau_*$}
\put(157,20){$\mathbf{R}$}
\end{picture}
\caption{Trace cohomology.}
\end{figure}
\begin{defn}\label{dfn2}
By an $i$-th trace cohomology  group $H^i_{tr}(V)$  of   $V$   one 
understands the  abelian subgroup  of   $\mathbf{R}$ defined by map (\ref{eq7}).
\end{defn}
\begin{rem}\label{rmk4}
The abelian group  $H^i_{tr}(V)$  is  called  a {\it pseudo-lattice},   see  
[Manin 2004]  \cite{Man1},  Section 1.    
The  endomorphisms  in  the category of pseudo-lattices are given by
multiplication of its points  by the real numbers $\alpha$
such that $\alpha H^i_{tr}(V)\subseteq H^i_{tr}(V)$.    It is known that the ring
$End~(H^i_{tr}(V))\cong \mathbf{Z}$ or
\linebreak
 $End~(H^i_{tr}(V))\otimes \mathbf{Q}$
is a   real algebraic number field.    In the latter case  $H^i_{tr}(V)\subset  End~(H^i_{tr}(V))\otimes \mathbf{Q}$,
see  [Manin 2004]  \cite{Man1}, Lemma 1.1.1 for the case of quadratic fields.  
Notice that  one can write multiplication  by $\alpha$ in a matrix form by fixing a basis in the pseudo-lattice;  
thus  the ring $End~(H^i_{tr}(V))$   is  a  commutative subring of  the matrix ring $M_{b_i}(\mathbf{Z})$,  
 where $b_i$ is equal to the  rank  of  pseudo-lattice,  i.e.  the cardinality of its basis.   
\end{rem}
\begin{rem}\label{rmk5}
 Notice that the trace cohomology $H^i_{tr}(V)$ is an abelian group with order,
see  [Goodearl 1986]  \cite{G}  for an introduction. The total order is defined 
by an order-preserving homomorphism $H^i_{tr}(V)\to \mathbf{R}$ given by formula
(\ref{eq7}). 
 \end{rem}

\subsection{Weil's Conjectures}
Let ${\Bbb F}_q$ be a finite field with $q=p^r$ elements and $V:=V({\Bbb F}_q)$ be a
smooth $n$-dimensional projective variety over ${\Bbb F}_q$.  The famous  Weil conjectures  
establish a deep relation between the arithmetic of  $V$ and topology of
the variety $V_\mathbf{C}$ defined by the polynomial equations over the field of complex numbers
 [Weil 1949]  \cite{Wei1}.    Namely,  let  $N_m$ be the number of
rational points of $V$ over the field ${\Bbb F}_{q^m}$ and 
\begin{equation}
Z_V(t)=\exp~\left(\sum_{m=1}^{\infty} N_m {t^m\over m}\right)
\end{equation}
the corresponding zeta function.   Weil conjectured  that:  (i) $Z_V(t)$ is a quotient 
of polynomials with rational coefficients; (ii) $Z_V(q^{-n}t^{-1})=\pm q^{n{\chi\over 2}} t^{\chi} Z_V(t)$,
where $\chi$ is the Euler-Poincar\'e characteristic of $V_\mathbf{C}$;  (iii) $Z_V(t)$ satisfies  an analog
of the Riemann Hypothesis,   i.e.    
\begin{equation}\label{eq2}
Z_V(t)={P_1(t) P_3(t)\dots P_{2n-1}(t)\over P_0(t) P_2(t)\dots P_{2n}(t)},
\end{equation}
so that $P_0(t)=1-t^r, P_{2n}(t)=1-q^n t$ and for each $1\le i\le 2n-1$ 
the polynomial $P_i(t)$ has integer coefficients and can be written 
in the form $P_i(t)=\prod (1-\alpha_{ij}t)$,  where $\alpha_{ij}$
are algebraic integers with  $|\alpha_{ij}|=q^{i\over 2}$;   (iv)  the degree
of polynomial $P_i(t)$ is equal to the $i$-th  Betti number of variety 
$V_\mathbf{C}$.  
The properties (i)-(iv)  are true for algebraic curves (i.e. for  $n=1$)
and Weil pointed out that  they  follow from a 
cohomology theory of the  variety $V$.  
Such a cohomology was constructed by Grothendieck and called  the  $\ell$-adic
cohomology;  all but conjecture (iii)  can be deduced from basic properties
of the $\ell$-adic cohomology   [Grothendieck 1968]  \cite{Gro1}.

\subsection{Kuga-Sato varieties}
Let $\Gamma(N)$ be the principal congruence subgroup of level $N\ge 3$. 
Denote by $X(N)=\mathbb{H}/\Gamma(N)$ the corresponding modular curve,
where $\mathbb{H}:=\{z=x+iy\in \mathbf{C} ~|~ y>0\}$ is the Lobachevsky half-plane. 
The  {\it Kuga-Sato variety} of level $N$   is the $k$-th power of the universal 
elliptic curve $\mathcal{E}$  over the modular curve, i.e. 
\begin{equation}
V_N=\underbrace{\mathcal{E}\times_{X(N)}\dots\times _{X(N)}
\mathcal{E}}_{k ~\hbox{times}}.
\end{equation}
In what follows,  we assume that the variety $V_N$ is compact.
Such a compactification is described in   [Deligne 1969]  \cite[Lemma 5.4]{Del1}.  
It is known, that
\begin{equation}
H^i_{et}(V_N; {\bf Q}_{\ell})\cong S_{i+1}(\Gamma(N)),
\end{equation}
where  $H^i_{et}(V_N; {\bf Q}_{\ell})$ is the $\ell$-adic cohomology 
of $V_N$ and $S_{i+1}(\Gamma(N))$ is the space of cusp forms for 
the group $\Gamma(N)$ [Deligne 1969]  \cite[Definition 2.8 and Theorem 2.10]{Del1}  and [Scholl 1985]  \cite[Section 2.5]{Sch1}.

\section{Proof of theorem \ref{thm1}}
For the sake of clarity,  let us outline main ideas. 
The trace cohomology $H^i_{tr}(V)$ will be used to construct a positive-definite Hermitian
form $\varphi(x,y)$ on the cohomology group $H^i(V(k))$;  see also remark \ref{rmk5}.  Such a construction 
involves the Deligne-Scholl theory linking the $\ell$-adic cohomology $H^i_{et}(V; ${\bf Q}$_{\ell})$
of the Kuga-Sato variety $V$  with the space of cusp forms $S_{i+1}(\Gamma)$ of weight $i+1$ for a finite index subgroup
$\Gamma\subset SL_2(\mathbf{Z})$,  see [Deligne 1969]  \cite{Del1}  and [Scholl 1985]  \cite{Sch1}.
It is proved that the Petersson inner product on $S_{i+1}$ defines,  via the trace cohomology,  
the required form $\varphi(x,y)$.  Since the regular maps of $V$ preserve  the form 
$\varphi(x,y)$ modulo a positive constant,  one obtains  an analog of the Riemann hypothesis for the zeta function
of $V$.    We shall split the proof in a series of lemmas. 
\begin{lem}\label{lem1}
{\bf (Deligne-Scholl)}
If $V$ is the Kuga-Sato variety, 
then  there exists a finite index subgroup $\Gamma$
of  the modular group $SL_2(\mathbf{Z})$,  such that
\begin{equation}
\dim  H^i_{tr}(V)=2\dim_\mathbf{C}  S_{i+1}(\Gamma),
\end{equation}
where $S_{i+1}(\Gamma)$ is the  space of cusp forms of 
weight $i+1$  relatively  group $\Gamma$.
\end{lem}
\begin{proof}
This lemma follows from the results of  [Deligne 1969]  
\cite[Definition 2.8 and Theorem 2.10]{Del1}  and [Scholl 1985]  \cite[Section 2.5]{Sch1}.
Namely,   let $\Gamma$ be a finite index subgroup of $SL_2(\mathbf{Z})$,  such that the modular
curve $X_{\Gamma}:={\Bbb H}/\Gamma$ can be defined over the field $\mathbf{Q}$.
It was proved that for each prime $\ell$ there exists a continuous homomorphism
\begin{equation}
\rho: Gal~(\overline{\mathbf{Q}} ~|~ \mathbf{Q})\to End~(W),
\end{equation}
where $W$ is a $2d$-dimensional vector space over $\ell$-adic numbers {\bf Q}$_{\ell}$ and 
\linebreak
$d=\dim_\mathbf{C} S_{i+1}(\Gamma)$,  see [Scholl 1985]  \cite{Sch1}.
It was proved earlier,  that $W\cong H^i_{et}(V; ${\bf Q}$_{\ell})$
 for a variety $V$ over the field $\mathbf{Q}$ 
and  the (arithmetic) Frobenius element of the Galois group  $Gal~(\overline{\mathbf{Q}} ~|~ \mathbf{Q})$
corresponds to the (geometric) Frobenius  endomorphism of the $\ell$-adic cohomology
$H^i_{et}(V; ${\bf Q}$_{\ell})$,  see [Deligne 1969]  \cite{Del1} for $\Gamma$ being 
a congruence group.

Let $V(k)$ be a variety over the complex numbers associated to $V$.  By the comparison
theorem
\begin{equation}
H^i_{et}(V; \hbox{{\bf Q}}_{\ell})\otimes_{\hbox{{\bf Q}}_{\ell}} \mathbf{C}\cong H^i(V(k); \mathbf{C}),
\end{equation}
see e.g. [Hartshorne 1977]  \cite[p. 454]{H}.   On the other hand, from definition \ref{dfn2}
we have   $\dim  ~H_{tr}^i(V)=\dim ~H^i(V(k);  \mathbf{C})$  and,  therefore,
\begin{equation}
\dim~H_{tr}^i(V)=\dim~H_{et}^i(V).   
\end{equation}
But $\dim  H^i_{et}(V; ${\bf Q}$_{\ell})=2\dim_\mathbf{C} S_{i+1}(\Gamma)$  by the 
Deligne-Scholl theory   and,  therefore,
\begin{equation}\label{eq11}
\dim  H^i_{tr}(V)=2\dim_\mathbf{C} S_{i+1}(\Gamma).
\end{equation}
Lemma \ref{lem1} follows.  
\end{proof}


\begin{lem}\label{lem2}
The trace cohomology $H^i_{tr}(V)$ defines a  $\mathbf{Z}$-module
embedding 
\begin{equation}
H^i(V(k))\hookrightarrow S_{i+1}(\Gamma). 
\end{equation}
\end{lem}
\begin{proof}
Recall that the {\it Petersson inner product} 
\begin{equation}
(x,y) ~: ~S_{i+1}(\Gamma)\times  S_{i+1}(\Gamma)\to \mathbf{C}
\end{equation}
on the space $S_{i+1}(\Gamma)$  is given by the integral
\begin{equation}
(f,g)=\int_{X_{\Gamma}} f(z)\overline{g(z)} (\Im z)^{i+1} dz. 
\end{equation}
The product   is linear in $f$ and conjugate-linear in $g$,  so that
 $(g,f)=\overline{(f,g)}$  and $(f,f)>0$ for all $f\ne 0$,
 see e.g.   [Milne 1997]  \cite[p. 57]{MI}.  

Fix a basis $\{\alpha_1,\dots,\alpha_d; \beta_1,\dots, \beta_d\}$ in  the 
$\mathbf{Z}$-module $H_{tr}^i(V)$.  In view of the standard properties
of scalar product $(x,y)$,   there exists a unique cusp form $g\in S_{i+1}(\Gamma)$,
such that 
\begin{equation}\label{eq16}
(f_j, g)=\alpha_j+i\beta_j,
\end{equation}
where $\{f_1,\dots, f_d\}$ is the orthonormal basis in  $S_{i+1}(\Gamma)$ 
consisting of the Hecke eigenforms.  
But $\alpha_j$ and $\beta_j$ are the image of generators of the $\mathbf{Z}$-module
$H^i(V(k))$ under the trace map $\tau_*$,  see definition \ref{dfn2}.   Using formula
 (\ref{eq16}),  one   defines an embedding 
\begin{equation}
\iota: ~H^i(V(k))\hookrightarrow S_{i+1}(\Gamma),
\end{equation}
whose image $\iota(H^i(V(k)))$ is a $\mathbf{Z}$-module generated by the 
real and imaginary parts of the Hecke  eigenforms $f_j\in S_{i+1}(\Gamma)$.
Lemma \ref{lem2} follows.
\end{proof}

\begin{cor}\label{cor1}
There exists a unique positive-definite Hermitian form
\begin{equation}\label{eq18}
\varphi(x,y) ~: ~H^i(V(k))\times  H^i(V(k))\to \mathbf{C},
\end{equation}
on the $\mathbf{Z}$-module $H^i(V(k))$ coming from the Petersson inner product
on the space $S_{i+1}(\Gamma)$.  
\end{cor}
\begin{proof}
The Petersson inner product is a Hermitian form  because $(g,f)=\overline{(f,g)}$
and a positive-definite form  because  $(f,f)>0$ for all  $f\ne 0$.
It is easy to see, that such a form is unique. 
In view of lemma \ref{lem2},  one gets the conclusion of 
corollary \ref{cor1}. 
\end{proof}

\begin{lem}\label{lmk6}
The ring $End~(H_{tr}^i(V))$ is isomorphic to the ring ${\Bbb T}_{i+1}(\Gamma)$
of the Hecke operators on the space $S_{i+1}(\Gamma)$ and it is a commutative 
subring of the matrix ring  $End~(H^i(V(k)))$. 
\end{lem}
\begin{proof}
In view of lemma \ref{lem1} and remark \ref{rmk4},   the ring 
 $End~(H_{tr}^i(V))$ is generated by the eigenvalues of Hecke operators 
 corresponding to their (common) eigenform $f_j\in S_{i+1}(\Gamma)$; 
notice that   ${\Bbb T}_{i+1}(\Gamma)$  is always a non-trivial ring   if $\Gamma$
is a congruence subgroup and extends to such for the non-congruence subgroups
of finite index as shown in [Scholl 1985]  \cite{Sch1}.  
On the other hand,  it is known that the Hecke ring ${\Bbb T}_{i+1}(\Gamma)$ 
is isomorphic to a commutative subring of the matrix ring $End~(H^i(V(k)))$
represented by the symmetric matrices with positive integer entries, see e.g. 
[Milne 1997]  \cite{MI}.   Lemma \ref{lmk6}  follows. 
\end{proof}

\begin{lem}\label{lem3}
Each  regular map $f: V\to V$ induces a linear map 
 $f^i_*: H^i (V(k))\to H^i (V(k))$ of degree $\deg(f^i_*)$,
 whose characteristic polynomial $char (f^i_*)$ has integer 
 coefficients and roots of the absolute value  $|\lambda|=[\deg(f^i_*)]^{1\over 2n}$. 
\end{lem}
\begin{proof}
Consider a regular map $f: V\to V$ obtained by the reduction modulo $q$ of 
an algebraic map $\widetilde f: V(k)\to V(k)$ of the corresponding variety
over the field of complex numbers. 
Such a map always exists, see  [Hartshorne 2010] \cite[Theorem 22.1]{H2}.
Let us show that the linear map 
$f^i_*: H^i(V(k))\to H^i(V(k))$ induced by $\widetilde f$ on the integral 
cohomology $H^i(V(k))$ must preserve, up to  a constant multiple, the 
positive-definite Hermitian form $\varphi(x,y)$ on $H^i(V(k))$ given by 
formula (\ref{eq18}).     

Indeed,  let $\widetilde f(V(k))\subseteq V(k)$ be a constructible subset; 
clearly,  such a subset  carries  the structure of an algebraic variety. 
We  repeat the trace cohomology construction for the variety $\widetilde f(V(k))$;
thus one gets a positive-definite Hermitian form $\widetilde\varphi (x,y)$ on
$H^i(\widetilde f(V(k)))$.   But $H^i(\widetilde f(V(k)))\subseteq H^i(V(k))$
and therefore one gets yet another positive-definite Hermitian form $\varphi(x,y)$
on $H^i(\widetilde f(V(k)))$.  Since such a form is unique (see lemma \ref{lem2}),
one concludes that $\widetilde\varphi(x,y)$ coincides with $\varphi(x,y)$ modulo
a positive factor $C$.  It is easy to see,  that $C=[\deg(f^i_*)]^{1\over n}$.
Indeed,  the volume form can be calculated by the formula $v=|\det~(f^i_*)|v_0=
\deg(f^i_*)v_0$;  on the other hand,  the  multiplication by $C$ map gives the
volume $v=C^nv_0$, where $n$ is the dimension of variety $V$.

Let $\lambda$ be a root of the characteristic polynomial $char (f^i_*):=\det~(\lambda I-f^i_*)$.
Since the kernel of the map $\lambda I- f^i_*$ is non-trivial,  let $x\in Ker~(\lambda I-f^i_*)$ 
be a non-zero element;  clearly,  $f^i_* ~x=\lambda x$.  Consider the value of scalar 
product $(x,y)=\varphi(x,y)$ on $x=y=f^i_* x$,  i.e.
\begin{equation}\label{eq19}
(f^i_* ~x, f^i_* ~x)=(\lambda x,\lambda x)=\lambda\bar\lambda (x,x).
\end{equation}
On the other hand,
\begin{equation}\label{eq20}
(f^i_* ~x, f^i_* ~x)=[\deg(f^i_*)]^{1\over n} (x,x).
\end{equation}
Because $(x,x)\ne 0$,  one can cancel it in (\ref{eq19}) and (\ref{eq20}), so that
\begin{equation}\label{eq21}
\lambda\bar\lambda=[\deg(f^i_*)]^{1\over n}  \quad\hbox{or} \quad |\lambda|=[\deg(f^i_*)]^{1\over 2n}.
\end{equation}
Note that $char (f^i_*)\in \mathbf{Z}[\lambda]$ because $H^i(V(k))$ is a $\mathbf{Z}$-module;
lemma \ref{lem3} follows. 
\end{proof}

\begin{rem}\label{rmk7}
Note that any non-trivial map 
$f^i_*\in End~(H^i_{tr}(V))$ 
\linebreak
$\subset End~(H^i(V(k)))$
corresponds to a non-algebraic  (transcendental) map $\widetilde f: V(k)\to V(k)$,
because the roots of $char (f^i_*)$ are real numbers in this case.  Of course,
there are many other examples of the non-algebraic maps  $\widetilde f: V(k)\to V(k)$.
\end{rem}

\begin{lem}\label{lem4}
$\deg(f^i_*)=[\deg(f)]^i$.
\end{lem}
\begin{proof}
It is well known,  that the cusp forms $g(z)\in S_{i+1}(\Gamma)$ are bijective with the
holomorphic differentials
\begin{equation}\label{eq22}
g(z)dz^{i+1\over 2}
\end{equation}
on the Riemann surface $X_{\Gamma}={\Bbb H}/\Gamma$.   To prove lemma \ref{lem4},
one can use the Riemann-Hurwitz formula:
\begin{equation}\label{eq23}
2g(Y)-2=m  ~[2g(X)-2]+\sum_{P} (e_P-1),
\end{equation}
where $e_P$ is the multiplicity at the point $P$ of an $m$-fold holomorphic map
$Y\to X$ between the Riemann surfaces of genus $g(Y)$ and $g(X)$,  see 
e.g.   [Milne 1997]  \cite[p. 17]{MI}.  Because the differential (\ref{eq22}) is locally
defined,  one can substitute in (\ref{eq23}) $g(X)=g(Y)=0$  and assume $P=0$
to be a unique ramification point.  Thus 
\begin{equation}\label{eq24}
m={e_P+1\over 2}
\end{equation}
and the $m$-fold differential (\ref{eq22}) implies $e_P=i$,  i.e. the holomorphic 
map $Y\to X$ is given by the formula 
\begin{equation}\label{eq25}
z\longmapsto z^i.
\end{equation}
On the other hand,  for a regular  map  $f: V\to V$ it holds $\deg(f)=\deg(\widetilde f)=
\deg(f^1_*)$.   Since degree is a multiplicative function on composition of maps,
one gets from (\ref{eq25}) and the  link between $i$-th cohomology of $V(k)$ and 
the space $S_{i+1}(\Gamma)$,  that 
\begin{equation}\label{eq26}
\deg(f^i_*)=[\deg(f^1_*)]^i=[\deg(f)]^i.
\end{equation}
Lemma \ref{lem4} follows.
\end{proof}

\begin{cor}\label{cor2}
{\bf (Riemann Hypothesis)}
The roots $\alpha_{ij}$ of polynomials $P_i(t)$ in formula (\ref{eq2})  
are algebraic numbers of the absolute value  $|\alpha_{ij}|=q^{{i\over 2}}$. 
\end{cor}
\begin{proof}
It is easy to see, that the Frobenius map  $f: ~(z_1,\dots, z_n)\mapsto (z_1^q,\dots, z_n^q)$
of variety $V$ is regular and $\deg(f)=q^n$.   Therefore,  one can apply lemmas \ref{lem3} 
and \ref{lem4} to  such a map and get the equality  $|\alpha_{ij}|=q^{{i\over 2}}$ for each 
$0\le i\le 2n-1$.   Corollary \ref{cor2} follows. 
\end{proof}

\bigskip
Corollary \ref{cor2}  finishes the proof of theorem \ref{thm1}.

\section{Examples}
The  groups $H^i_{tr}(V)$  are  truly   concrete and simple;
in this section we calculate   the trace cohomology for  $n=1$,  i.e.  when $V$
is a smooth  algebraic curve.  In particular,  we find the cardinality of  the set $\mathcal{E}({\Bbb F}_{q})$ 
 obtained by the reduction modulo $q$   of   an elliptic curve  with  complex multiplication.
\begin{ex}\label{exm1}
The trace cohomology of smooth  algebraic curve $\mathcal{C}({\Bbb F}_q)$ of
genus $g\ge 1$ is given by the formulas: 
\begin{equation}\label{eq32}
\left\{
\begin{array}{lll}
H_{tr}^0(\mathcal{C}) &\cong&  \mathbf{Z},\\
&&\\
H_{tr}^1(\mathcal{C}) &\cong&  \mathbf{Z}+\mathbf{Z}\theta_1+\dots+\mathbf{Z}\theta_{2g-1},\\
&&\\
H_{tr}^2(\mathcal{C}) &\cong&  \mathbf{Z},   
\end{array}
\right.
\end{equation}
where $\theta_i\in \mathbf{R}$ are  algebraically independent  integers
of a number field of degree $2g$.  
\end{ex}
\begin{proof}
It is known that the Serre $C^*$-algebra of the (generic) complex algebraic curve $\mathcal{C}$
is isomorphic to a  {\it toric} $AF$-algebra  ${\Bbb A}_{\theta}$,  
see \cite{Nik3} for the notation and details.   Moreover,  up to  a scaling constant $\mu>0$,
it holds
\begin{equation}\label{eq33}
\tau_*(K_0({\Bbb A}_{\theta}\otimes \mathcal{K}))=
\begin{cases}
\mathbf{Z}+\mathbf{Z}\theta_1, &  \hbox{if}  ~g=1\cr
             \mathbf{Z}+\mathbf{Z}\theta_1+\dots+\mathbf{Z}\theta_{6g-7}, & \hbox{if} ~g>1,
             \end{cases}
\end{equation}
where constants $\theta_i\in \mathbf{R}$ parametrize the moduli (Teichm\"uller) space of curves $\mathcal{C}$
 \cite{Nik3}.
If $\mathcal{C}$ is defined over a number field $k$,  then each $\theta_i$ is algebraic and 
their total number   is equal to  $2g-1$.  (Indeed,  since $Gal~(\bar k ~|~ k)$ acts on the torsion
points of $\mathcal{C}(k)$,  it is easy to see that the endomorphism ring of $\mathcal{C}(k)$ is non-trivial.  
Because such a ring  is isomorphic to the endomorphism ring of the Jacobian $Jac~\mathcal{C}$ 
and $\dim_\mathbf{C} Jac~\mathcal{C}=g$,  one concludes that $End ~\mathcal{C}(k)$ is a $\mathbf{Z}$-module
of rank $2g$ and each $\theta_i$ is an algebraic number.)  After scaling by a constant $\mu>0$,   one gets   
\begin{equation}\label{eq34}
H^1_{tr}(\mathcal{C}):=\tau_*(K_0({\Bbb A}_{\theta}\otimes \mathcal{K}))=
\mathbf{Z}+\mathbf{Z}\theta_1+\dots+\mathbf{Z}\theta_{2g-1}
\end{equation}
Because $H^0(\mathcal{C})\cong H^2(\mathcal{C})\cong \mathbf{Z}$,
one obtains   the rest of formulas (\ref{eq32}).
\end{proof}

\begin{rem}\label{rmk10}
If $k\cong \mathbf{Q}$,  then $\Gamma\cong\Gamma(N)$ is the principal congruence
subgroup of level $N$,   since $\mathcal{C}(\mathbf{Q})\cong X_{\Gamma(N)}$ for some integer $N$.
As explained,  the Petersson inner product on $S_2(\Gamma(N))$ gives rise to a positive-definite
Hermitian form $\varphi$ on the cohomology group $H^1(\mathcal{C})\cong \mathbf{Z}^{2g}$.
Note that the  form $\varphi$ can be obtained from the classical Riemann's bilinear relations
for the periods of curve $\mathcal{C}$;   this  yields Weil's proof of the Riemann hypothesis 
for function $Z_\mathcal{C}(t)$.   
\end{rem}
\begin{rem}
Notice   the cardinality of the set $\mathcal{C}({\Bbb F}_q)$ is  given  to  the formula
\begin{equation}\label{eqn30}
|\mathcal{C}({\Bbb F}_{q})|=  1+q-tr~(\omega)= 1+q-\sum_{i=1}^{2g}\lambda_i, 
\end{equation}
where $\lambda_i$ are  the  eigenvalues of the Frobenius endomorphism 
$\omega\in End~(H^1_{tr}(\mathcal{C}))$.  
\end{rem}
\begin{ex}\label{exm2}
The case $g=1$ is particularly  instructive;  for the sake of clarity,  we shall consider
elliptic curves having complex multiplication.  Let $\mathcal{E}({\Bbb F}_q)$  
be the reduction modulo $q$ of an elliptic with complex multiplication 
by the ring of integers of an imaginary quadratic field $\mathbf{Q}(\sqrt{-d})$,
see e.g.  [Silverman 1994]   \cite{S},  Chapter 2.  It is known,  that in this case 
the trace cohomology formulas (\ref{eq32})  take the form  
\begin{equation}\label{eq38}
\left\{
\begin{array}{lll}
H_{tr}^0(\mathcal{E}({\Bbb F}_q)) &\cong&  \mathbf{Z},\\
&&\\
H_{tr}^1(\mathcal{E}({\Bbb F}_q)) &\cong&  \mathbf{Z}+\mathbf{Z}\sqrt{d},\\
&&\\
H_{tr}^2(\mathcal{E}({\Bbb F}_q)) &\cong&  \mathbf{Z}.     
\end{array}
\right.
\end{equation}
We shall denote by  $\psi(\mathfrak{P})\in \mathbf{Q}(\sqrt{-d})$  the Gr\"ossencharacter
of  the prime ideal $\mathfrak{P}$ over $p$,   see  [Silverman 1994]   \cite{S},   p. 174.
It is easy to see, that in this case  the Frobenius endomorphism $\omega\in End~(H_{tr}^1(\mathcal{E}({\Bbb F}_q)))$
is given by the formula
\begin{equation}\label{eq39}
\omega={1\over 2}\left[\psi(\mathfrak{P})+\overline{\psi(\mathfrak{P})}\right] + 
{1\over 2}\sqrt{\left(\psi(\mathfrak{P})+\overline{\psi(\mathfrak{P})}\right)^2+4q}   
\end{equation}
and the corresponding eigenvalues 
\begin{equation}\label{eq40}
\left\{
\begin{array}{lll}
\lambda_1 &=& \omega= {1\over 2}\left[\psi(\mathfrak{P})+\overline{\psi(\mathfrak{P})}\right] + 
{1\over 2}\sqrt{\left(\psi(\mathfrak{P})+\overline{\psi(\mathfrak{P})}\right)^2+4q},  \\
\lambda_2 &=& \bar\omega= {1\over 2}\left[\psi(\mathfrak{P})+\overline{\psi(\mathfrak{P})}\right] -
{1\over 2}\sqrt{\left(\psi(\mathfrak{P})+\overline{\psi(\mathfrak{P})}\right)^2+4q}  .
\end{array}
\right.
\end{equation}
Using formula (\ref{eqn30}),  one  gets the following equation 
\begin{equation}\label{eqn41}
|\mathcal{E}({\Bbb F}_{q})|=  1-(\lambda_1+\lambda_2)+q=1-\psi(\mathfrak{P})-\overline{\psi(\mathfrak{P})} +q,
\end{equation}
which coincides with the well-known expression for  $|\mathcal{E}({\Bbb F}_{q})|$ in terms 
of the Gr\"ossencharacter,  see e.g.  [Silverman 1994]   \cite[p. 175]{S}. 
\end{ex}


\subsection*{Acknowledgment}
 The author would like to thank the anonymous referee who provided useful and detailed comments on a earlier version of the manuscript.

\end{document}